\documentclass[aap,MSNbibl,dvips]{arximspdf}

\doi{10.1214/12-AAP918} %kopijuoti is PTS
\volume{24}
\issue{1}
\pubyear{2014}
\firstpage{131}
\lastpage{149}

\makeatletter
\newtheorem{theorem}{Theorem}
\newtheorem{lemma}[theorem]{Lemma}
\newproclaim{remark}{Remark}
\makeatother

\begin{document}
\begin{frontmatter}

\title{Minimising MCMC variance via diffusion limits, with~an
application to simulated tempering}
\runtitle{Minimising MCMC variance via diffusion limits}

\begin{aug}
\author[A]{\fnms{Gareth O.} \snm{Roberts}\ead[label=e1]{g.o.roberts@lancaster.ac.uk}}
\and
\author[B]{\fnms{Jeffrey S.} \snm{Rosenthal}\corref{}\ead[label=e2]{jeff@math.toronto.edu}\ead[label=u1,url]{http://probability.ca/jeff/}\thanksref{t1}}
\runauthor{G. O. Roberts and J. S. Rosenthal}
\affiliation{University of Warwick and University of Toronto}
\address[A]{Department of Statistics\\
University of Warwick\\
Coventry, CV4 7AL\\
United Kingdom\\
\printead{e1}} %adresu isvedimo komanda gale!
\address[B]{Department of Statistics\\
University of Toronto\\
Toronto, Ontario, M5S 3G3\\
Canada\\
\printead{e2}\\
\printead{u1}}
\end{aug}
\thankstext{t1}{Supported in part by NSERC of Canada.}

% HISTORY:
\received{\smonth{3} \syear{2012}}
\revised{\smonth{12} \syear{2012}}

% ABSTRACT
%
\begin{abstract}
We derive new results comparing the asymptotic variance of diffusions
by writing them as appropriate limits of discrete-time birth--death
chains which themselves satisfy Peskun orderings. We then apply our
results to simulated tempering algorithms to establish which choice of
inverse temperatures minimises the asymptotic variance of all
functionals and thus leads to the most efficient MCMC algorithm.
\end{abstract}

% KEYWORDS

\end{frontmatter}

%s1 #&#
\section{Introduction}\label{sec-intro}

Markov chain Monte Carlo (MCMC) algorithms are very widely used to
approximately compute expectations with respect to complicated
high-dimensional distributions; see, for example,
\cite{tierney,mcmchandbook}. Specifically, if a Markov chain $\{X_n\}$
has stationary distribution $\pi$ on state space $\mathcal{X}$, and
$h\dvtx \mathcal{X}\to\mathbf{R}$ with $\pi|h|<\infty$, then $\pi(h):=
\int h(x)  \pi(dx)$ can be estimated by $\frac{1}{n} \sum_{i=1}^n
h(X_i)$ for suitably large $n$. This estimator is unbiased if the chain
is started in stationarity (i.e., if $X_0 \sim\pi$), and in any case
has bias only of order $1/n$. Furthermore, it is consistent provided
the Markov chain is $\phi$-irreducible. Thus, the efficiency of the
estimator is often measured in terms of the asymptotic variance
$\operatorname{Var}_\pi(h,P):= \lim_{n\to\infty} \frac{1}{n}
\operatorname{Var}_\pi ( \sum_{i=1}^n h(X_i)  )$ (where the
subscript $\pi$ indicates that $\{X_n\}$ is in stationarity): the
smaller the variance, the better the estimator.

An important question in MCMC research is how to \textit{optimise} it, that is,
how to choose the Markov chain optimally; see, for example,
\cite{geyerstatsci,mira}. This leads to the question of how to
\textit{compare} different Markov chains. Indeed, for two different
\mbox{$\phi$-}irreducible Markov chain kernels $P_1$ and $P_2$ on
$\mathcal{X}$, both having the same invariant probability measure
$\pi$, we say that $P_1$ \textit{dominates} $P_2$ \textit{in the efficiency ordering},
written $P_1 \succeq P_2$, if $\operatorname{Var}_\pi(h,P_1)
\leq\operatorname{Var}_\pi(h,P_2)$ for all $L^2(\pi)$ functionals
$h\dvtx \mathcal{X}\to\mathbf{R}$, that is, if $P_1$ is ``better'' than
$P_2$ in the sense of being uniformly more efficient for estimating
expectations of functionals.

It was proved by Peskun~\cite{peskun} for finite state spaces, and by
Tierney~\cite{tierney2} for general state spaces
(see also \cite{mirageyer,mira}), that if $P_1$ and $P_2$ are discrete-time
Markov chains which are both reversible with respect to the same
stationary distribution $\pi$, then a sufficient condition for $P_1
\succeq P_2$ is that $P_1(x,A) \geq P_2(x,A)$ for all $x\in\mathcal{X}$
and $A\in\mathcal{F}$ with $x\notin A$, that is, that $P_1$
\textit{dominates} $P_2$ \textit{off the diagonal}.

Meanwhile, diffusion limits have become a common way to establish
asymptotic comparisons of MCMC algorithms
\cite{RGG,lang,statsci,bedard1,bedard2,bedardrosenthal,stuart}.
Specifically, if $P_{1,d}$ and $P_{2,d}$ are two different Markov
kernels in dimension $d$ (for $d=1,2,3,\ldots$), with diffusion limits
$P_{1,*}$ and $P_{2,*}$ respectively as $d\to\infty$, then one way to
show that $P_{1,d}$ is more efficient than $P_{2,d}$ for large $d$ is
to prove that $P_{1,*}$ is more efficient that $P_{2,*}$. This leads to
the question of how to establish that one diffusion is more efficient
than another. In some cases
(e.g., random-walk Metropolis \cite{RGG}, and Langevin algorithms \cite{lang}), this is easy since one diffusion
is simply a time-change of the other. But more general diffusion
comparisons are less clear; for example, the processes' spectral gaps
\[
1 - \sup \biggl\{ \int h(y) P(x,dy) \dvtx \int h(y) \pi(dy) = 0, \int
h^2(y) \pi(dy) = 1 \biggr\}
\]
can be ordered directly by using Dirichlet forms, but this does not
lead to bounds on the asymptotic variances.

In this paper, we develop (Section~\ref{sec-compare}) a new comparison
of asymptotic variance of diffusions. Specifically, we prove
(Theorem~\ref{mainthm}) that if $P_i$ are Langevin diffusions with
respect to the same stationary distribution $\pi$, with variance
functions $\sigma_i^2$ (for $i=1,2$), then if $\sigma_1^2(x)
\geq\sigma_2^2(x)$ for all $x$, then $P_1 \succeq P_2$, that is, $P_1$
is more efficient than $P_2$. (We note that Mira and
Leisen~\cite{miral1,miral2} extended the Peskun ordering in an
interesting way to continuous-time Markov processes on finite state
spaces, and on general state spaces when the processes have generators
which can be represented as $G_i f(x) = \int f(y)  Q_i(x,dy)$ and which
satisfy the condition that $Q_1(x, A \setminus \{x\}) \geq Q_2(x, A
\setminus\{x\})$ for all $x$ and $A$. However, their results do not
appear to apply in our context, since generators of diffusions involve
differentiation and thus do not admit such representation.)

We then consider (Section~\ref{sec-simtemp}) simulated tempering
algorithms \cite{marinari,geyerstatsci}, and in particular the question
of how best to choose the intermediate temperatures. It was previously
shown in~\cite{yves}, generalising some results in the physics
literature~\cite{kofke1,PP}, that a particular choice of temperatures
(which leads to an asymptotic temperature-swap acceptance rate of
0.234) maximises the asymptotic $L^2$ jumping distance, that is,
$\lim_{n\to\infty} \mathbf{E}(|X_n - X_{n-1}|^2)$. (Indeed, this result
has already influenced adaptive MCMC algorithms for simulated
tempering; see, for example,~\cite{fort2011}.) However, the previous
papers did {not} prove a diffusion limit, nor did they provide any
comparisons of Markov chain variances. In this paper, we establish
(Theorem~\ref{simtempthm}) diffusion limits for certain simulated
tempering algorithms. We then apply our diffusion comparison results to
prove (Theorem~\ref{simtempoptthm}) that the given choice of
temperatures does indeed minimise the asymptotic variance of all
functionals.

%s2 #&#
\section{Comparison of diffusions}\label{sec-compare}

Let $\pi\dvtx \mathcal{X}\to(0,\infty)$ be a $C^1$ target density
function, where $\mathcal{X}$ is either $\mathbf{R}$ or some finite
interval $[a,b]$. We shall consider nonexplosive Langevin diffusions
$X^{\sigma}$ on $\mathcal{X}$ with stationary density $\pi$, satisfying
%
%
%e1 #&#
\begin{equation}
\label{first} dX^{\sigma}_t = \sigma \bigl(X^{\sigma}
\bigr)\,dB_t + \bigl( \tfrac{1}{2} \sigma^2
\bigl(X^{\sigma}_t \bigr) \log\pi'
\bigl(X^{\sigma}_t \bigr) +\sigma \bigl(X^{\sigma}_t
\bigr) \sigma' \bigl(X^{\sigma}_t \bigr) \bigr)\,dt
\end{equation}
for some $C^1$ function $\sigma\dvtx  \mathcal{X}\to[\underline{k},
\overline{k}]$ for some fixed $0<\underline{k} < \overline{k} <
\infty$, and with reflecting boundaries at $a$ and $b$ in the case
$\mathcal{X}=[a,b]$.

For two such diffusions $X^{\sigma_1}$ and $X^{\sigma_2}$, we write
(similarly to the above) that $X^{\sigma_1} \succeq X^{\sigma_2}$, and
say that $X^{\sigma_1}$ \textit{dominates} $X^{\sigma_2}$ \textit{in the efficiency
ordering}, if for all $L^2(\pi)$ functionals $f\dvtx
\mathcal{X}\to\mathbf{R}$,
\[
\lim_{T\to\infty} T^{-1/2} \operatorname{Var} \biggl( \int
_{0}^T f \bigl(X^{\sigma_1}_s
\bigr)\,ds \biggr) \leq \lim_{T\to\infty} T^{-1/2}
\operatorname{Var} \biggl( \int_{0}^T f
\bigl(X^{\sigma_2}_s \bigr)\,ds \biggr).
\]
%
%Because of ergodicity of the diffusions,
%this covariance ordering is independent of starting values [RIGHT?].

We wish to argue that if $\sigma_1(x) \geq\sigma_2(x)$ for all $x$,
then $X^{\sigma_1} \succeq X^{\sigma_2}$. Intuitively, this is because
$X^{\sigma_1}$ ``moves faster'' than $X^{\sigma_2}$, while maintaining
the same stationary distribution. Indeed, if $\sigma_1$ and $\sigma_2$
are constants, then this result is trivial
(and implicit in earlier works \cite{RGG,lang,statsci}), since then $X^{\sigma_1}_t$ has the
same distribution as $X^{\sigma_2}_{ct}$ where $c=\sigma_1/\sigma_2>1$;
that is, $X^{\sigma_1}$ accomplishes the same sampling as
$X^{\sigma_2}$ in a shorter time, so it must be more efficient.
However, if $\sigma_1$ and $\sigma_2$ are nonconstant functions, then
the comparison of $X^{\sigma_1}$ and $X^{\sigma_2}$ is less clear.

To make theoretical progress, we assume:

\begin{longlist}[(A1)]
\item[(A1)] $\pi$ is log-Lipschitz function on $\mathcal{X}$; that is,
    there is $L<\infty $ with
%
%e2 #&#
\begin{equation}\label{Ldef}
\bigl|\log\pi(y)-\log\pi(x)\bigr| \leq L |y-x|, \qquad x,y\in\mathcal{X}.
\end{equation}

\item[(A2)] Either (a) $\mathcal{X}$ is
    a bounded interval $[a,b]$, and the diffusions $X^\sigma$ have
    reflecting boundaries at $a$ and $b$, or (b) $\mathcal{X}$ is
    all of $\mathbf{R}$, and $\pi$ has exponentially-bounded tails;
    that is, there is $0<K<\infty$ and $r>0$ such that
\[
\pi(x+y) \leq \pi(x) e^{-ry}, \qquad x > K, y>0
\]
and
\[
\pi(x-y) \leq \pi(x) e^{-ry}, \qquad x < -K, y>0.
\]

In case (A2)(b), we can then find sufficiently large $q \geq K$ such
that
%
%e3 #&#
\begin{equation}\label{qdef}
\mathop{\sum_{i}}_{|i/m| \geq q}
\pi(i/m) \leq (1/4) \sum _i \pi(i/m)\qquad \mbox{for all } m\in
\mathbf{N} %%%%%%%%%%%%%%%\atop
\end{equation}
[where the sums in (\ref{qdef}) must be finite due to~(\ref{Ldef})],
and then set
%
%
%e4 #&#
\begin{equation}\label{Qdef}
Q = \inf \bigl\{ \pi(x) \dvtx |x| \leq q+1\bigr\},
\end{equation}
which must be positive by continuity of $\pi$ and compactness of the
interval $[-q-1, q+1]$.

Our main result is then the following.
\end{longlist}

%th1 #&#
\begin{theorem}\label{mainthm}
If $X^{\sigma_1}$ and $X^{\sigma_2}$ are two Langevin diffusions of the
form~(\ref{first}) with respect to the same density $\pi$, with
variance functions $\sigma_1$ and $\sigma_2$ respectively, and if
$\sigma_1(x) \geq\sigma_2(x)$ for all $x \in \mathcal{X}$, then
assuming \textup{(A1)} and \textup{(A2)}, we have $X^{\sigma_1}\succeq
X^{\sigma_2}$.
\end{theorem}

%s2.1 #&#
\subsection{\texorpdfstring{Proof of Theorem~\protect\ref{mainthm}}
{Proof of Theorem 1}}\label{sec-mainproof}

To prove Theorem~\ref{mainthm}, we introduce auxiliary processes for
each $m\in\mathbf{N}$. Given $\sigma\dvtx \mathcal{X}\to\mathbf{R}$,
let $S = 2 \overline{k} e^L$, and let $Z^{m,\sigma}$ be a discrete-time
birth and death process on the discrete state space $\mathcal{X}_m:=
\{i/m; i \in\mathbf{Z}\}$ in case~(A2)(b), or $\mathcal{X}_m:= \{i/m;
i \in\mathbf{Z}\} \cap[a,b]$ in case (A2)(a), with transition
probabilities given by
\begin{eqnarray*}
P \bigl(i/m, (i+1)/m \bigr) &=& {1 \over2S} \biggl( \sigma^2(i/m)
+ {\sigma^2((i+1)/m) \pi((i+1)/m) \over\pi(i/m)} \biggr),
\\
P \bigl(i/m, (i-1)/m \bigr) &=& {1 \over2S} \biggl(
\sigma^2(i/m) + {\sigma^2((i-1)/m) \pi((i-1)/m) \over\pi(i/m)} \biggr)
\end{eqnarray*}
and
\[
P(i/m, i/m) = 1 - P \bigl(i/m, (i+1)/m \bigr) - P \bigl(i/m, (i-1)/m \bigr).
\]
(In case (A2)(a), any transitions which would cause the process to move
out of the interval $[a,b]$ are instead given probability~0.) These
transition rates are chosen to satisfy detailed balance with respect to
the stationary distribution $\pi_m$ on $\mathcal{X}_m$ given by
$\pi_m(i/m) = \pi(i/m)/\sum_{x\in\mathcal{X}_m} \pi(x)$ [and $S$ is
chosen to be large enough to ensure that $P(i/m, (i+1)/m) + P(i/m,
(i-1)/m) \leq1$].

In terms of $Z^{m,\sigma}$, we then let $\{Y^\sigma_{m,t}\}_{t \geq0}$
be the continuous-time version of $Z^{m,\sigma}$, speeded up by a
factor of $m^2S/2$, that is, defined by $Y^\sigma_{m,t} =
Z^{m,\sigma}_{\lfloor m^2 S t / 2 \rfloor}$ for \mbox{$t \geq0$}. (Here
and throughout, $\lfloor r \rfloor$ is the floor function which rounds
$r$ down to the next integer, e.g. $\lfloor6.8 \rfloor= 6$ and
$\lfloor-2.1 \rfloor= -3$.) It then follows that $Y_{m,t}$ converges
to~$X^{m,\sigma}$, as stated in the following lemma (whose proof is
deferred until the end of the paper, since it uses similar ideas to
those of the following section).
%easily from standard arguments (see e.g. \cite{ethierkurtz,bhat})

%le2 #&#
\begin{lemma}\label{weakconvlemma}
Assuming \textup{(A1)} and \textup{(A2)}, as $m \rightarrow \infty$,
the processes $Y^\sigma_m$ converge weakly (in the Skorokhod topology)
to $X^{\sigma}$.
\end{lemma}

We then apply the usual discrete-time Peskun ordering to the
$Z^{m,\sigma}$ processes, as follows.

%le3 #&#
\begin{lemma}\label{orderlemma}%
Suppose that $\sigma_1(x) \geq\sigma_2(x)$ for all $x\in\mathbf{R}$.
Then $Z^{m,\sigma_1}\succeq Z^{m,\sigma_2}$.
\end{lemma}

\begin{pf}
By inspection, the fact that $\sigma_1(x) \geq\sigma_2(x)$ implies that
\[
\mathbf{P} \bigl(Z^{m,\sigma_1}_{(i+1)/m} = j+1 \mid Z^{m,\sigma_1}_{i/m}=j
\bigr) \geq\mathbf{P} \bigl(Z^{m,\sigma_2}_{(i+1)/m} = j+1 \mid Z^{m,\sigma
_2}_{i/m}=j \bigr)
\]
and
\[
\mathbf{P} \bigl(Z^{m,\sigma_1}_{(i+1)/m} = j-1 \mid Z^{m,\sigma_1}_{i/m}=j
\bigr) \geq \mathbf{P} \bigl(Z^{m,\sigma_2}_{(i+1)/m} = j-1 \mid
Z^{m,\sigma
_2}_{i/m}=j \bigr).
\]
It follows that $Z^{m,\sigma_1}$ dominates $Z^{m,\sigma_2}$ off the
diagonal. The usual discrete-time Peskun ordering
\cite{peskun,tierney2} thus implies that $Z^{m,\sigma_1}\succeq
Z^{m,\sigma_2}$.
\end{pf}

To continue, let
\[
V_*(f,\sigma):= \lim_{T\to\infty} T^{-1}
\operatorname{Var}_\pi \biggl( \int_0^T
f \bigl(X^{\sigma}_s \bigr)\,ds \biggr),
\]
which we assume satisfies the usual relation
\[
V_*(f,\sigma) = \int_{-\infty}^\infty
\operatorname{Cov}_\pi \bigl(f \bigl(X^{\sigma}_0 \bigr),f
\bigl(X^{\sigma}_s \bigr) \bigr)\,ds.
\]
Also, let
\[
V_m(f,\sigma):= \lim_{n\to\infty} n^{-1}
\operatorname{Var}_\pi \Biggl( \sum_{i=1}^{mn}
f \bigl(Z^{m,\sigma
}_i \bigr) \Biggr),
\]
which we assume satisfies the usual relation
\[
V_m(f,\sigma) = \sum_{i=-\infty}^\infty
\operatorname{Cov}_\pi \bigl(f \bigl(Z^{m,\sigma
}_0 \bigr),f
\bigl(Z^{m,\sigma}_i \bigr) \bigr).
\]
(In both cases, the subscript $\pi$ indicates that the process is
assumed to be in stationarity, all the way from time $-\infty$ to
$\infty$.) We then have the following.

%le4 #&#
\begin{lemma}\label{Vconvlemma}%
% Let $\pi:\IR\to(0,\infty)$ be continuous and integrable.
% Let $\pi_m(i/m) = \pi(i/m)/\sum_{j\in\IZ} \pi(j/m)$.
% Suppose for each $m \in\IN$, $X_m:= \{X_{m,n}\}_{n=0}^\infty$ is a
% discrete-time reversible process on the state space
% $\X_m:= \{i/m : i\in\IZ\}$, with stationary distribution $\pi_m$.
% Let $Y_m$ be a continuous-time version of $X_m$, defined by
% $Y_{m,t} = X_{\lfloor mt \rfloor}$ for $t \geq0$.
% Suppose the processes $Y_m$ converge weakly (in the Skorokhod
% topology) to a diffusion process $Z$ having stationary density
% proportional to $\pi$.
Let $G_m$ be the spectral gap of the process $Z^{m,\sigma}$. Assume
there is some constant $g>0$ such that $G_m \geq g/m^2$ for all $m$.
Then for all bounded functions $f:\mathbf{R}\to\mathbf{R}$,
$\lim_{m\to\infty} (m^2S/2) V_m(f,\sigma) = V_*(f,\sigma)$.
\end{lemma}

\begin{pf}
%Since it's continuous time, the spectrum and auto-covariances are
%nonnegative.
Let
\[
A_{m,t} = \operatorname{Cov}_\pi \bigl[f
\bigl(Z^{m,\sigma}_0 \bigr),f \bigl(Z^{m,\sigma
}_{\lfloor m^2St/2 \rfloor} \bigr) \bigr]
\]
and let
\[
A_{*,t} = \operatorname{Cov}_\pi \bigl[f
\bigl(X^{\sigma}_0 \bigr),f \bigl(X^{\sigma}_t
\bigr) \bigr].
\]
Then
\[
V_*(f,\sigma) = \int_{\infty}^\infty A_{*,t}\,dt
\]
and (since $\lfloor m^2St/2 \rfloor$ is a step-function of $t$, with
steps of size $m^2 S/2$)
\[
V_m(f,\sigma) = { \int_{-\infty}^\infty A_{m,t}\,dt
\over m^2 S/2 }.
\]
Now, by Lemma~\ref{weakconvlemma}, since $f$ is bounded,
\[
\lim_{m\to\infty} A_{m,t} = A_{*,t}.
\]
To continue, let $F$ be the forward operator corresponding to the chain
$Z^{m,\sigma}$, that is, $F h(x) = \mathbf{E}[h(Z^{m,\sigma}_1)\mid
Z^{m,\sigma}_0=x]$. Then since $F$ is reversible, it follows from
Lemma~2.3 of~\cite{liu-bka94} that
\[
\bigl\|F^t\bigr\| = \|F\|^t = \sup \bigl\{ \operatorname{Cov}_\pi
\bigl[h_1 \bigl(Z^{m,\sigma}_0 \bigr),
h_2 \bigl(Z^{m,\sigma}_t \bigr) \bigr] \dvtx
\operatorname{Var}_\pi(h_1) = \operatorname{Var}_\pi(h_2)
= 1 \bigr\}.
\]
Letting $v = \operatorname{Var}_\pi[f(X)]$, we then have, for all
$m\in\mathbf {N}$ and $t \geq0$, that
\begin{eqnarray*}
A_{m,t} &=& \operatorname{Cov}_\pi \bigl[f \bigl(Z^{m,\sigma}_0
\bigr),f \bigl(Z^{m,\sigma }_{\lfloor m^2St/2 \rfloor}\bigr) \bigr]
\\
&\leq& \sup \bigl\{ \operatorname{Cov}_\pi \bigl[h
\bigl(Z^{m,\sigma}_0 \bigr),h \bigl(Z^{m,\sigma}_{\lfloor
m^2St/2 \rfloor}
\bigr) \bigr] \dvtx h\in L^2(\pi), \operatorname{Var}_\pi
\bigl[h(X) \bigr]=v \bigr\}
\\
&=& v \big\|F^{\lfloor m^2St/2 \rfloor}\big\| = v \|F\|^{\lfloor m^2St/2 \rfloor} = v (1-G_m)^{\lfloor m^2St/2 \rfloor}
\\
&\leq& v \bigl(1-g/m^2 \bigr)^{\lfloor m^2St/2 \rfloor} \leq v
\bigl(e^{-g/m^2} \bigr)^{m^2St/2} = v e^{-gSt/2}.
\end{eqnarray*}
Hence,
\[
V_m(f,\sigma) = \int_{-\infty}^\infty
A_{m,t}\,dt \leq 2 \int_0^\infty
A_{m,t}\,dt % EDITED SLIGHTLY BY JSR RIGHT AFTER SUBMITTING TO AAP.
\leq 4v/gS < \infty.
\]
Hence, by the dominated convergence theorem,
\[
\lim_{m\to\infty} \int_{-\infty}^\infty
A_{m,t}\,dt = \lim_{m\to\infty} \int_{-\infty}^\infty
A_{*,t}\,dt,
\]
that is,
\[
\lim_{m\to\infty} \bigl(m^2S/2 \bigr) V_m(f,
\sigma) = V_*(f,\sigma)
\]
as claimed.\vadjust{\goodbreak}
\end{pf}

To make use of Lemma~\ref{Vconvlemma}, we need to bound the spectral
gaps of the $Z^{m,\sigma}$ processes. We do this using a capacitance
argument; see, for example,~\cite{sinclair}. Let
\[
\kappa_m = \mathop{\inf_{A\subseteq\mathcal{X}_m}}_{0<\pi(A) \leq1/2} {1 \over\pi_m(A)}
\sum_{x\in A} P_m \bigl(x,A^C\bigr) \pi_m(x)
\]
be the capacitance of $Z^{m,\sigma}$. We prove

%le5 #&#
\begin{lemma}\label{gaplemma}%
The capacitance $\kappa_m$ satisfies that
\[
\kappa_m \geq \min \biggl( {k e^{-L} r \over2m},
{Qk e^{-2L/m} \over2m} \biggr),
\]
where the quantities $L$ and $Q$ are defined in~(\ref{Ldef})
and~(\ref{Qdef}), respectively, and where the bound reduces to simply
$\kappa_m \geq{k e^{-L} r \over2m}$ in case \textup{(A2)(a)}.
\end{lemma}

\begin{pf}
We consider two different cases [only the second of which can occur in
case (A2)(a)]:

\begin{longlist}[(i)]
\item[(i)] $\exists a \in A$ with $|a| \leq q$. Then,
    since $\pi_m(A) \leq 1/2$, there is $j\in\mathbf{Z}$ with $|j/m|
    \leq q$ and $j/m \in A$ and either $(j+1)/m \in A^C$ or $(j-1)/m
    \in A^C$. Assume WOLOG that $(j+1)/m \in A^C$. We will need the
    following estimate on $\sum_{j \in\mathbf{Z}} \pi(j/m)$. For $x
    \in[i/m, (i+1)/m)$,
\[
\pi(x) \geq \pi(i/m) e^{-L(x - i/m) }
\]
so that
\begin{eqnarray*}
\int_{i/m}^{(i+1)/m} \pi(x) &\geq& \pi(i/m) \int _0^{1/m} e^{-Lu}\,du
= \pi(i/m) \biggl( {1 - e^{-L/m} \over L} \biggr)
\\
&=& \pi(i/m) e^{-L/m} \biggl( {e^{L/m} - 1 \over L} \biggr)
\geq \pi(i/m) e^{-L/m} \biggl({L/m \over L}\biggr)
\\
&=& {\pi(i/m)   e^{-L/m} \over m}.
\end{eqnarray*}
Therefore summing both sides over all $i \in\mathbf{Z}$,
\[
1 = \int_{-\infty}^\infty\pi(x)\,dx \geq
{e^{-L/m} \over m} \sum_{i \in\mathbf{Z}} \pi(i/m),
\]
whence
\[
\sum_{i \in\mathbf{Z}} \pi(i/m) \leq m e^{L/m}.
\]
Then
\begin{eqnarray*}
\sum_{x\in A} P_m \bigl(x,A^C \bigr) \pi_m(x) &\geq& \pi_m(j/m)
P_m \bigl(j/m,(j+1)/m \bigr)
\\
&=& \pi_m(j/m) (1/2) \sigma^2(j/m) e^{-L/m}
\\
&\geq& \bigl(\pi(j/m)/m \bigr) (k/2) e^{-2 L/m}
\\
&\geq& Qk e^{- 2L/m}/2m.
\end{eqnarray*}
%
%[NEED TO FIX UP $\pi_m$ NORMALISING CONSTANT! SHOULD BE OKAY, IT'S
%BASICALLY A RIEMANN SUM APPROXIMATION. USE OF ``$L$'' SHOULD HELP.]

\item[(ii)] $A \subseteq(-\infty,q) \cup(q,\infty)$.
    Let $a\in A$ with $\pi(a) = \max\{\pi(x) \dvtx  x \in A\}$. Assume
    WOLOG that $a>0$. Then
\begin{eqnarray*}
\sum_{x\in A} P_m \bigl(x,A^C
\bigr) \pi_m(x) &\geq& \pi_m(a) P_m\bigl(a,a-(1/m) \bigr)
\\
&\geq& k e^{-L/m} \pi(a) \Big/ \mathop{\sum_{i}}_{|i/m| \geq a} \pi(i/m)
\\
&\geq& k e^{-L/m} \pi(a) \Big/ \Biggl[2 \sum_{j=0}^\infty
\pi(a) e^{-rj/m} \Biggr]
\\
&=& \frac{1}{2}k e^{-L/m} \bigl[1-e^{-r/m} \bigr] \leq
\frac{1}{2}k e^{-L} (r/m).
\end{eqnarray*}
\end{longlist}

Thus, in either case, the conclusion of the lemma is satisfied.
\end{pf}

Now, it is known (e.g.,~\cite{sinclair}) that the spectral gap can be
bounded in terms of the capacitance, specifically that $G_m
\geq\kappa_m^2/2$. Thus, for $m \geq1$,
\begin{eqnarray*}
G_m &\geq& \bigl[\min \bigl(\tfrac{1}{2}k e^{-L} (r/m), Qk
e^{-2L/m}/2m \bigr) \bigr]^2/2
\\
&\geq& \bigl[\min \bigl(\tfrac{1}{2}k e^{-L} (r/m), Qk
e^{-2L}/2m \bigr) \bigr]^2/2
\\
&=& g/m^2,
\end{eqnarray*}
where $g = [\min(\frac{1}{2}k e^{-L} r, Qk e^{-2L}/2)]^2/2 > 0$. This
together with Lemma~\ref{weakconvlemma} shows that the conditions of
Lemma~\ref{Vconvlemma} are satisfied. Hence, by Lemma~\ref{Vconvlemma},
$\lim_{m\to\infty} (m^2S/2) V_m(f,\sigma) = V_*(f,\sigma)$ for all
bounded functions $f$.

On the other hand, by Lemma~\ref{orderlemma}, $Z^{m,\sigma_1}\succeq
Z^{m,\sigma_2}$, that is, $V_m(f,\sigma_1) \leq V_m(f,\sigma_2)$.
Hence, for all bounded functions $f$,
%
%
%e5 #&#
\begin{eqnarray}\label{boundedf}
V_*(f,\sigma_1) &=& \lim_{m\to\infty} \bigl(m^2S/2 \bigr)
V_m(f,\sigma_1)\nonumber
\\
&\leq& \lim_{m\to\infty} \bigl(m^2S/2 \bigr) V_m(f,
\sigma_2)
\\
 &=& V_*(f,\sigma_2). \nonumber
\end{eqnarray}

Finally, if $f$ is in $L^2$ but not bounded, then letting
\[
f_m(x) = \cases{m, &\quad $f(x) > m$,
\cr
f(x), &\quad $-m \leq f(x)\leq m$,
\cr
-m, &\quad $f(x) < -m$,}
\]
we have by the monotone (or dominated) convergence theorem that
$V_*(f,\break \sigma_1) = \lim_{m\to\infty} V_*(f_m,\sigma_1)$ and
$V_*(f,\sigma_2) = \lim_{m\to\infty} V_*(f_m,\sigma_2)$. Hence, it
follows from~(\ref{boundedf}) that $V_*(f,\sigma_1) \leq
V_*(f,\sigma_2)$ for all $L^2(\pi)$ functions $f$. That is,
$X^{\sigma_1} \succeq X^{\sigma_2}$, thus proving
Theorem~\ref{mainthm}.

%s3 #&#
\section{Simulated tempering diffusion limit}\label{sec-simtemp}

We now apply our results to a version of the simulated tempering
algorithm. Specifically, following~\cite{yves}, we consider a
$d$-dimensional target density
%
%
%e6 #&#
\begin{equation}\label{prodform}
f_d(x) = e^{d K} \prod_{i=1}^d f(x_i)
\end{equation}
for some unnormalised one-dimensional density function $f\dvtx
\mathbf{R}\to[0,\infty)$, where $K = - \log(\int f(x)\,dx)$ is the
corresponding normalising constant. (Although~(\ref{prodform}) is a
very restrictive assumption, it is
known~\cite{RGG,statsci,bedard1,bedard2,bedardrosenthal,stuart} that
conclusions drawn from this special case are often approximately
applicable in much broader contexts.) We consider simulated tempering
in $d$ dimensions, with inverse-temperatures chosen as follows:
$\beta_0^{(d)}=1$, and $\beta_{i+1}^{(d)} = \beta_i^{(d)} -
{\ell(\beta_i^{(d)}) \over d^{1/2}}$ for some fixed $C^1$ function
$\ell\dvtx  [0,1] \to\mathbf{R}$. (The question then becomes, what is
the optimal choice of $\ell$.) As for when to stop adding new
temperature values, we fix some $\chi\in(0,1)$ and keep going until the
temperatures drop\vspace*{-2pt} below $\chi$; that is, we stop at
temperature $\beta_{k(d)}^{(d)}$ where $k(d) = \sup\{i \dvtx
\beta_i^{(d)} \geq\chi\}$.

We shall consider a joint process $(y_n^{(d)},X_n)$, with $X_n
\in\mathbf{R}^d$, and with $y_n^{(d)} \in E_d:= \{\beta^{(d)}_i; 0\leq
i \leq k(d)\}$ defined as follows. If $y_{n-1}=\beta_i^{(d)}$ [where
$0<i<k(d)$], then the chain proceeds by choosing $X_{n-1} \sim
f^\beta$, then proposing $Z_n$ to be $\beta_{i+1}$ or $\beta_{i-1}$
with probability $1/2$ each, and finally accepting $Z_n$ with the usual
Metropolis acceptance probability. (A proposed move to
$\beta_{-1}^{(d)}$ or $\beta_{k(d)+1}^{(d)}$ is automatically
rejected.) We assume, as in~\cite{yves}, that the chain then
immediately jumps to stationary at the new temperature, that is, that
mixing within a temperature is infinitely more efficient than mixing
between temperatures.

The process $(y_n^{(d)}, X_n)$ is thus a Markov chain on the
state space $ E_d \times\mathbf{R}^d$, with joint stationary density
given by
\[
f_d(\beta,x) = e^{d K(\beta)} \prod_{i=1}^d
f^\beta(x_i),
\]
where $K(\beta) = - \log\int f^\beta(x)\,dx$ is the normalising
constant.\vadjust{\goodbreak}

We now prove that the $\{y_n^{(d)}\}$ process has a diffusion limit
(similar to random-walk Metropolis and Langevin algorithms, see \cite{RGG,lang,statsci}), and furthermore the asymptotic variance of
the algorithm is minimised by choosing the function $\ell$ that leads
to an asymptotic temperature acceptance rate $\doteq0.234$.
Specifically, we prove the following:

%th6 #&#
\begin{theorem}\label{simtempthm}%
Under the above assumptions, the $\{y_n^{(d)}\}$ inverse-tempera\-ture
process, when speeded up by a factor of $d$, converges in the Skorokhod
topology as $d\to\infty$ to a diffusion limit $\{X_t\}_{t \geq0}$
satisfying
%
%
%e7 #&#
\begin{eqnarray}\label{diffusiondef}
dX_t &=& \biggl[ 2 \ell^2 \Phi \biggl( {- \ell I^{1/2} \over2}
\biggr) \biggr] ^{1/2}\,dB_t
\nonumber\\[-8pt]\\[-8pt]
&&{} + \biggl[ \ell(X) \ell' (X) \Phi \biggl(
{-I^{1/2} \ell\over2} \biggr) - \ell^2 \biggl({ \ell I^{1/2} \over2}
\biggr)'
\phi \biggl( {-I^{1/2} \ell\over2} \biggr) \biggr]\,dt\nonumber % + LOCAL TIME TERMS!!!
\end{eqnarray}
for $X_t$ in $(\chi, 1)$ with reflecting boundaries at both $\chi$ and
$1$. Furthermore, the speed of this diffusion is maximised, and the
asymptotic variance of all $L^2$ functionals is minimised, when the
function $\ell$ is chosen so that the asymptotic temperature acceptance
rate is equal to~0.234 (to three decimal places).
\end{theorem}

Then, combining Theorems~\ref{mainthm}~and~\ref{simtempthm},
we immediately obtain: %[RIGHT?]:

%th7 #&#
\begin{theorem}\label{simtempoptthm}%
For the above simulated tempering algorithm, for any $L^2$ functional $f$,
the choice of $\ell$ which minimises the limiting
asymptotic variance $V_*(f) = \lim_{m\to\infty} V_m(f)$, is the same
as the choice which maximises $\sigma(x)$, that is, is the
choice which leads to an asymptotic temperature acceptance probability
of 0.234 (to three decimal places).
\end{theorem}

\begin{remark*}
In this context, it was proved in~\cite{yves} that as $d\to\infty$, the
choice of~$\ell$ leading to an asymptotic temperature acceptance rate
$\doteq0.234$ maximises the expected squared jumping distance of the
$\{y_n^{(d)}\}$ process. However, the question of whether that choice
would also minimise the asymptotic variance for any $L^2$ function was
left open. That question is resolved by Theorem~\ref{simtempoptthm}.
\end{remark*}

%s3.1 #&#
\subsection{\texorpdfstring{Proof of Theorem~\protect\ref{simtempthm}}
{Proof of Theorem 6}}

The key computation for proving Theorem~\ref{simtempthm} will be given
next, but first we require some additional notation. We let
$\mathrm{int}(E_d)$ denote $E_d\setminus\{1, \beta^{(d)}_{k(d)}\}$. We
also denote by $G^{(d)}$ the generator of the inverse-temperature
process $\{y_n^{(d)}\}$ and set $H$ to be the set of all functions $h
\in C^2[\chi, 1]$ with $h'(\chi) = h'(1)=0$. We also let $G^*$ be the
generator of the diffusion given in (\ref{diffusiondef}), defined, for
all functions $h \in H$, by
%
%
%e8 #&#
\begin{equation}\label{Geqn}
G^*h = { \sigma^2(x) h''(x) \over2} + \mu(x) h'(x), \qquad h \in H,
\end{equation}
where
\[
\mu(x) = \ell(x) \ell' (x) \Phi \biggl( {-I^{1/2} \ell\over2}
\biggr) - \ell^2 \biggl({ \ell I^{1/2} \over2} \biggr)'
\phi \biggl( {-I^{1/2} \ell\over2} \biggr)
\]
and
%
%e9 #&#
\begin{equation}\label{sigmaeqn}
\sigma^2(x) = 2 \ell^2 \Phi \biggl(
{- \ell I^{1/2} \over2} \biggr).
\end{equation}

To proceed, we apply the powerful weak convergence theory of
\cite{ethierkurtz}. We do this using a technique for limiting
reflecting processes similar to the arguments in Ward and
Glynn~\cite{glynn}. We first note that by page~17 and Chapter~8 of
\cite{ethierkurtz}, the set $\{(h,G^*h); h \in H\}$ forms a core for
the generator of the diffusion process described above in
(\ref{diffusiondef}) (i.e., the closure of the restriction of the
generator to that set is again equal to the generator itself). Hence,
by Theorems 1.6.1~and~4.2.11 of \cite{ethierkurtz}, we need to show
that, for any pair $(h, G^*h)$ with $h \in H$, there exists a sequence
$(h_d, dG^{(d)}h_d)_{d\in\mathbf{N}}$ such that
%
%
%e10 #&#
\begin{equation}\label{unifconv}
\lim_{d\to\infty} \sup_{x \in E_d} \bigl|h(x) - h_d(x)\bigr| =
0
\end{equation}
and
%
%e11 #&#
\begin{equation}\label{needtoprove}
\lim_{d\to\infty} \sup_{x \in E_d} \bigl|G^*h(x) - d
G^{(d)}h_d(x)\bigr| = 0.
\end{equation}

To establish this convergence on $\mathrm{int}(E_d)$, we can simply let
$h_d=h$ (see Lemma~\ref{simtemplemma} below). However, to establish the
convergence on the boundary of $E_d$
(Lemma~\ref{simtempboundarylemma}), we need to modify $h$ slightly
[without destroying the convergence on $\mathrm{int}(E_d)$]. We do this
as follows. First, given any $h \in H$, we let
\[
\overline{h}_d (x) = h \bigl( \gamma_d(x) \bigr),
\]
where
\[
\gamma_d(x) = { (1 - \chi) x + \chi- \chi_d \over1 - \chi_d },
\]
so that $\overline{h}_d$ is just like $h$ except ``stretched'' to be
defined on $[\chi_d,   1]$ instead of just on $[\chi,   1]$.
% We shall be particularly interested in setting $a=\chi_d$, which has
%the
% effect of mapping $E_d$ precisely into $[\chi, 1]$.
Here we set $\chi_d=\beta_{k(d)}^{(d)}$, and
$\chi_d^+=\beta_{k(d)-1}^{(d)}$; thus $\chi_d \leq\chi\leq\chi_d^+$.
Notice that since $\chi_d \rightarrow\chi$ as $d \to\infty$,
$\overline{h}_d$ and its first and second derivatives converge to~$h$
and its corresponding derivatives uniformly for $x \in[\chi_d, 1]$ as
$d \rightarrow\infty$.

Finally, given the function $h$, we let $\eta(x)$ to be any smooth
function: $[\chi, 1] \to\mathbf{R}$ satisfying
\[
\eta'(\chi) = h'' (\chi)\ell(\chi)/2\quad\mbox{and}\quad \eta'(1) = h''(1) \ell(1)/2
\]
and then set
\[
h_d (x) = \overline{h}_d(x) + d^{-1/2} \eta
\bigl(\gamma_d(x) \bigr) = h \bigl(\gamma_d(x) \bigr) +
d^{-1/2} \eta \bigl(\gamma_d(x) \bigr),
\]
so that $h_d(x)$ is similar to $\overline{h}_d(x)$ except with the
addition of a separate $O(d^{-1/2})$ term (which will only be relevant
at the boundary points, i.e., in Lemma~\ref{simtempboundarylemma}
below). In particular, (\ref{unifconv}) certainly holds.

In light of the above discussion, Theorem~\ref{simtempthm} will follow
by establishing~(\ref{needtoprove}), which is done in
Lemmas~\ref{simtemplemma} and~\ref{simtempboundarylemma} below.

%le8 #&#
\begin{lemma}\label{simtemplemma}%
For all $h \in H$,
%
%
%e12 #&#
\begin{equation}\label{simtemph}
\lim_{d\to\infty} \sup_{x \in\mathrm{int}(E_d)} \bigl|d G^{(d)} h(x)
- G^* h(x) \bigr| = 0
\end{equation}
and
%
%e13 #&#
\begin{equation}\label{simtemphd}
\lim_{d\to\infty} \sup_{x \in\mathrm{int}(E_d)} \bigl|d
G^{(d)}h_d(x) - G^* h(x) \bigr| = 0.
\end{equation}
%
%for all bounded $C^2$ functions $h:[0,1] \to\IR$ with $h'(1)=h'(
%where $G^*$ is the generator for a diffusion satisfying~
%and furthermore the convergence is uniform.
\end{lemma}

\begin{pf}
%Let $h\dvtx [0,1] \to\IR$ be a bounded $C^2$ function with
%$h'(1)=h'(\chi)=0$.
We begin with a Taylor series expansion for $G^{(d)}$. Since the
computations shall get somewhat messy, we wish to keep only
higher-order terms, so for simplicity we shall use the notation $
\stackrel{r(d)}{\approx}$ to mean that the expansion holds up to
terms of order $1/r(d)$, uniformly for $x \in E_d$, as $d
\rightarrow\infty$ [e.g., $\mathit{LHS} \stackrel{d}{\approx}  \mathit{RHS}$ means
that $\lim_{d\to\infty} \sup_{x \in E_d} d(\mathit{LHS} - \mathit{RHS}) = 0$]. Then for
bounded $C^2$ functionals $h$, we have (combining the two $h''$ terms
together) that for $\beta_i^{(d)} \in \mathrm{int}(E_d)$:
\begin{eqnarray*}
G^{(d)}h \bigl(\beta_i^{(d)} \bigr) &\stackrel{d}{\approx}&
{h'(\beta_i^{(d)}) \over2} \bigl[\alpha^+ \bigl( \beta_{i+1}^{(d)} -
\beta_i^{(d)} \bigr) + \alpha^- \bigl(\beta _{i-1}^{(d)} -
\beta_i^{(d)} \bigr) \bigr]
\\[1pt]
&&{} + {h''(\beta_i^{(d)}) \over2} \bigl[
\bigl(\beta_{i+1}^{(d)} - \beta_i^{(d)} \bigr) ^2 \alpha^+ \bigr]
\\[1pt]
&\stackrel{d}{\approx}& {h'(\beta_i^{(d)}) \over2} \bigl[\alpha^+ \bigl(
\beta_{i+1}^{(d)} - \beta_i^{(d)} \bigr) +
\alpha^- \bigl(\beta _{i-1}^{(d)} - \beta_i^{(d)}
\bigr) \bigr]
\\[1pt]
&&{} + {h''(\beta_i^{(d)}) \over2} \bigl[ \bigl(\beta_{i+1}^{(d)}
- \beta_i^{(d)} \bigr) ^2 \alpha^+ \bigr]
\\[1pt]
&=& {h'(\beta_i^{(d)}) \over2}{\alpha^- \ell(\beta_{i-1}^{(d)} ) -\alpha^+
\ell(\beta_i^{(d)}) \over d^{1/2} }
\\[1pt]
&&{} + {h''(\beta_i^{(d)}) \over2}
\biggl[ {\ell(\beta_i^{(d)}) ^2 \alpha^+ \over d} \biggr],
\end{eqnarray*}
where $\alpha^+$ is the probability of accepting an upwards move, and
$\alpha^-$ is the probability of accepting a downwards move.\vadjust{\goodbreak}
%$$
%=: \mu(\beta_i^{(d)}) h'(\beta_i^{(d)} ) +
%$$
%(Here, and elsewhere, $\approx$ is understood as an expansion which
%holds
%uniformly for $x \in E_d$, as $d \rightarrow\infty$. Therefore in
%this case,
%$\sup_{x \in E_d} d(LHS - RHS ) \rightarrow0$ as $d \rightarrow

To continue, we let $g = \log f$, and
\[
M(\beta) = \mathbf{E}^\beta(g) = {\int\log f (x) f^\beta(x)\,dx \over\int f^\beta(x)\,dx}
\]
and
\[
I(\beta) = \operatorname{Var}^\beta(g) = {\int(\log f (x) )^2 f^\beta(x)\,dx \over
\int f^\beta(x)\,dx } - M(
\beta)^2.
\]
It follows, as in~\cite{yves}, that $M'(\beta) = I(\beta)$ and
$K'(\beta) = -M(\beta)$, so $K''(\beta)=-M'(\beta)=-I(\beta)$. We also
define $\overline{g} = g - M(\beta)$.

For shorthand, we write $\beta= \beta_i^{(d)}$, and $\ell=
\ell(\beta_{i}^{(d)})$, and ${\underline{\ell}}=
\ell(\beta_{i-1}^{(d)})$, and ${\underline{\varepsilon}}=
\beta_{i-1}^{(d)}-\beta_i^{(d)} = {\underline{\ell}}/d^{1/2}$, and
$\varepsilon= \beta_{i}^{(d)}-\beta_{i+1}^{(d)} = \ell/d^{1/2}$, and $I
= I(\beta)$ and $K'' = K''(\beta)$ and $K''' = K'''(\beta)$.

Then, with $X \sim f^\beta$,
%
%
%e14 #&#
\begin{eqnarray} \label{alphaminus}
\alpha^- &=& \mathbf{E} \biggl[ 1 \wedge{
f_d^{\beta+{\underline{\varepsilon}}}(X) e^{d K(\beta+{\underline
{\varepsilon}})} \over f_d^{\beta}(X) e^{d K(\beta)} } \biggr]\nonumber
\\
&=& \mathbf{E} \Biggl[ 1 \wedge \exp \Biggl( \bigl( K(\beta+{\underline{
\varepsilon}})-K(\beta) \bigr) d + {\underline{\varepsilon}}d M(\beta) + {
\underline{\varepsilon}} \sum_{i=1}^d
\overline{g}(X_i) \Biggr) \Biggr]
\nonumber\\
&\stackrel{d^{1/2}}{\approx}& \mathbf{E} \biggl[ 1 \wedge \exp \biggl(
{d {\underline{\varepsilon}}^2 \over2} K'' + {d {\underline{\varepsilon}}^3 \over6}
K''' + N \bigl( 0, I {\underline{\varepsilon}}^2 d \bigr) \biggr)
\biggr]\nonumber\\[-8pt]\\[-8pt]
&=& \mathbf{E} \biggl[ 1 \wedge \exp \biggl( {{\underline{\ell}}^2 \over2}
K'' + {{\underline{\varepsilon}}{\underline{\ell}}^2 \over6} K'''
+ N \bigl( 0, I {\underline{\ell}}^2 \bigr) \biggr) \biggr]\nonumber
\\
&=& \Phi \biggl( - {I^{1/2} {\underline{\ell}}\over2} + {{\underline{\varepsilon}}{\underline{\ell}}K''' \over6 I^{1/2}} \biggr)\nonumber
\\
&&{} + \exp \bigl({\underline{\varepsilon}} {\underline{\ell}}^2K''' / 6 \bigr) \Phi \biggl( -
{I^{1/2} {\underline{\ell}}\over2} - {{\underline{\varepsilon}}{\underline{\ell}}K''' \over6 I^{1/2}} \biggr).\nonumber
\end{eqnarray}

Similarly,
\begin{eqnarray*}
\alpha^+ &=& \mathbf{E} \biggl[ 1 \wedge{
f_d^{\beta-\varepsilon}(X) e^{d K(\beta-\varepsilon)} \over
f_d^{\beta}(X) e^{d K(\beta)} } \biggr]
\\
&=& \mathbf{E} \Biggl[ 1 \wedge \exp \Biggl( \bigl( K(\beta-\varepsilon)-K( \beta)
\bigr) d - \varepsilon d M(\beta) - \varepsilon\sum_{i=1}^d
\overline{g}(X_i) \Biggr) \Biggr]
\\
&\stackrel{1}{\approx}& \mathbf{E} \biggl[ 1 \wedge \exp \biggl(
{d \varepsilon^2 \over2} K'' - N \bigl( 0, I
\varepsilon^2 d \bigr) \biggr) \biggr]
\\
&=& \mathbf{E} \biggl[ 1 \wedge \exp \biggl( {\ell^2 \over2} I -
{\varepsilon\ell^2 \over6} K''' - N \bigl(
0, I \ell^2 \bigr) \biggr) \biggr]
\\
&=& \Phi \biggl( - {I^{1/2} \ell\over2} - {\varepsilon\ell K''' \over6 I^{1/2}} \biggr)
\\
&&{} + \exp \bigl(-\varepsilon\ell^2 K'''
/ 6 \bigr) \Phi \biggl( - {I^{1/2} \ell\over2} - {\varepsilon\ell K''' \over6 I^{1/2}}
\biggr).
\end{eqnarray*}
Hence
\begin{eqnarray*}
\alpha^+ \bigl(\beta_i^{(d)} \bigr) & \stackrel{d^{1/2}}{\approx}& \Phi \biggl(-{I^{1/2}(\beta_i^{(d)} ) \ell\over2} -
{\varepsilon \ell K'''(\beta_i^{(d)} ) \over6 I^{1/2}(\beta_i^{(d)} )}
\biggr)
\\
&&{}+ \exp \bigl(-\varepsilon\ell^2 \bigl(\beta_i^{(d)}
\bigr) K'''(\beta_i )/6
\bigr)\Phi \biggl(-{I^{1/2}(\beta_i^{(d)} ) \ell\over2} + {\varepsilon
\ell K'''(\beta_i^{(d)} ) \over6 I^{1/2}(\beta_i^{(d)} )}
\biggr).
\end{eqnarray*}
A first order approximation of this expression is
\[
\alpha^+ \bigl(\beta_i^{(d)} \bigr) \stackrel{1}{\approx}
2 \Phi \biggl(-{I^{1/2}(\beta_i^{(d)} ) \ell\over2} \biggr).
\]

Next, we note that in the current setting, $\beta$ is itself
marginally a Markov chain with uniform stationary distribution among
all temperatures. In fact it is a birth and death process, and hence
reversible. So, by detailed balance,
\[
\alpha^- = \alpha^+ \bigl(\beta_i^{(d)} - \ell/\sqrt d
\bigr).
\]
Therefore,
\begin{eqnarray*}
\alpha^- \bigl(\beta_i^{(d)} \bigr) &=& \alpha^+ \bigl(
\beta_i^{(d)} - \ell/\sqrt d \bigr)
\\
&\stackrel{d^{1/2}}{\approx}& \alpha^+ \bigl(\beta_i^{(d)}
\bigr)
\\
&&{} - {(\ell(\beta_i^{(d)}) I^{1/2} (\beta_i^{(d)}) )' \over2} \biggl( {-
\ell\over\sqrt d} \biggr) \phi
\biggl(-{I^{1/2}(\beta_i^{(d)} ) \ell\over2} - {\varepsilon
\ell K'''(\beta_i^{(d)} ) \over6 I^{1/2}(\beta_i^{(d)} )} \biggr)
\\
&&{}- \exp \bigl(-\varepsilon\ell^2 \bigl(
\beta_i^{(d)} \bigr) K'''(
\beta_i )/6 \bigr){(\ell(\beta_i^{(d)}) I^{1/2} (\beta_i^{(d)}) )' \over2}
\\
&&\hphantom{-}{}\times \biggl(
{- \ell\over\sqrt d} \biggr)\phi \biggl(-{I^{1/2}(\beta_i^{(d)} )
\ell\over2} +
{\varepsilon
\ell K'''(\beta_i^{(d)} ) \over6 I^{1/2}(\beta_i^{(d)} )} \biggr).
\end{eqnarray*}
Then, since ${\underline{\ell}} \stackrel{d^{1/2}}{\approx} \ell+
{\underline{\varepsilon}}\ell' \stackrel{d^{1/2}}{\approx} \ell+
\varepsilon\ell' = \ell+ {\ell\ell' \over d^{1/2}}$, we compute that
\begin{eqnarray*}
\mu \bigl(\beta_i^{(d)} \bigr) &\stackrel{d^{1/2}}{\approx}&
{1 \over2d^{1/2}} \biggl[ - \alpha^+ \ell + \biggl(\ell+ { \ell\ell' \over
d^{1/2} } \biggr)
\\
&&\hspace*{30pt}{} \times \biggl( \alpha^+ \bigl(\beta_i^{(d)} \bigr) \\
&&\hspace*{47.2pt}{}-
{(\ell(\beta_i^{(d)}) I^{1/2} (\beta_i^{(d)}) )' \over2}
\\
&&\hspace*{12.6pt}\hspace*{46pt}{}\times \biggl( {-
\ell\over\sqrt d} \biggr)\phi \biggl(-
{I^{1/2}(\beta_i^{(d)} ) \ell\over2} - {\varepsilon
\ell K'''(\beta_i^{(d)} ) \over6 I^{1/2}(\beta_i^{(d)} )} \biggr)
\\
&&\hspace*{12.6pt}\hspace*{46pt}{} - \exp \bigl(-\varepsilon\ell^2 \bigl(\beta_i^{(d)}
\bigr) K'''(\beta_i )/6
\bigr){(\ell
(\beta_i^{(d)}) I^{1/2} (\beta_i^{(d)}) )' \over2}
\\
&&\hspace*{12.6pt}\hspace*{73pt}{}\times \biggl( {- \ell
\over\sqrt d} \biggr)\phi
\biggl(-{I^{1/2}(\beta_i^{(d)} ) \ell
\over2} + {\varepsilon
\ell K'''(\beta_i^{(d)} ) \over6 I^{1/2}(\beta_i^{(d)} )} \biggr) \biggr) \biggr].
\end{eqnarray*}
Hence, ignoring all lower order terms,
\begin{eqnarray*}
\mu \bigl(\beta_i^{(d)} \bigr) &\stackrel{d^{1/2}}{\approx}&
{1 \over2 d^{1/2}} \biggl[ -\ell{(\ell(\beta_i^{(d)}) I^{1/2} (\beta_i^{(d)})
)' \over2}
\\
&&\hspace*{38pt}{}\times \biggl( {- \ell\over\sqrt d} \biggr) \phi
\biggl(-{I^{1/2}(\beta_i^{(d)} ) \ell\over2} - {\varepsilon \ell
K'''(\beta_i^{(d)} ) \over6 I^{1/2}(\beta_i^{(d)} )} \biggr)
\\
&&\hspace*{30pt}{}- \ell\exp{(\ell(\beta_i^{(d)}) I^{1/2} (\beta_i^{(d)}) )' \over 2}
\\
&&\hspace*{38pt}{}\times \biggl( {- \ell\over\sqrt d} \biggr) \phi
\biggl(-{I^{1/2}(\beta
_i^{(d)} ) \ell\over2} + {\varepsilon
\ell K'''(\beta_i^{(d)} ) \over6 I^{1/2}(\beta_i^{(d)} )} \biggr)
\\
&&\hspace*{112pt}{} + {2 \Phi (-{I^{1/2}(\beta_i^{(d)} ) \ell/2}  ) \ell\ell
' \over d^{1/2} } \biggr]
\\
&\stackrel{d^{1/2}}{\approx}& {1 \over d} \biggl[ -\ell^2
{(\ell(\beta_i^{(d)}) I^{1/2} (\beta_i^{(d)}) )' \over
2}\phi \biggl(-{I^{1/2}(\beta_i^{(d)} ) \ell\over2} \biggr)
\\
&&\hspace*{94Pt}{} + \Phi \biggl(-{I^{1/2}(\beta_i^{(d)} ) \ell\over2} \biggr) \ell\ell' \biggr].
\end{eqnarray*}

Similarly $\sigma^2 (\beta_i^{(d)})$ is to first order
\[
{2 \ell^2 \over d} \Phi \biggl(-{I^{1/2}(\beta_i^{(d)} ) \ell\over2} \biggr)
\]
so that we can write (for $0<\beta<1$)
\begin{eqnarray*}
G^dh &\stackrel{d}{\approx}& {1 \over d} \biggl( \ell^2 \Phi
\biggl(-{I^{1/2}(\beta_i^{(d)} ) \ell\over2} \biggr) h''(\beta)
\\
&&\hspace*{13pt}{}+ \biggl[ \Phi \biggl(-{I^{1/2}(\beta_i^{(d)} ) \ell\over2} \biggr) \ell\ell'
\\
&&\hspace*{30pt}{}  -\ell^2 {(\ell(\beta_i^{(d)}) I^{1/2} (\beta_i^{(d)})
)' \over2}\phi \biggl(-{I^{1/2}(\beta
_i^{(d)} ) \ell\over2} \biggr) \biggr]
 h'(\beta) \biggr).
\end{eqnarray*}
However, this expression is just $d^{-1}G^*h$, thus
establishing~(\ref{simtemph}).

Finally, to establish~(\ref{simtemphd}), we note that in this case the
terms $d^{-1/2} \eta(\gamma_d(x))$ and $\overline{h}_d(x) - h(x)$ are
both lower-order and do not affect the limit. Hence, (\ref{simtemphd})
follows directly from~(\ref{simtemph}).
% the statement about $\hchid(x)$ follows immediately.
% which for notational simplicity we now write simply as $h(x)$.
\end{pf}

%If this diffusion is then speeded up by a factor of $d$, then $Gh$
%is multiplied by $d$, i.e. the speeded-up diffusion's generator $G^*$
%satisfies
%$$
%G^*h \approx\ell^2 \Phim h''(\beta) +
% .
%$$
%$$ = G^* h $$
%But as is well-known, see for example Chapter 8 of \cite{ethierkurtz},
%a diffusion of the form $dX_t = \sigmai(X_t) B_t + \mu dt$
%has generator $Gh = (\sigma^2/2) h'' + \mu h'$.
%Hence, the generator $G^*$ is consistent with the diffusion
%$$ dX_t = \left[ 2 \ell^2 \Phi\left( {- \ell I^{1/2} \over2} \right)
%^{1/2} dB_t + \left[ \ell(X) \ell' (X) \Phi\left( {-I^{1/2} \ell
%{-I^{1/2} \ell\over2}

%It remains to clarify the uniformity of this convergence, and the
%convergence at the \textit{boundary}. We do this following Ward and
%Glynn~\cite{glynn}, as follows. [GARETH TO FILL THIS IN!]

The uniformity over $\mathrm{int}(E_d)$ for $h$ (as opposed to $h_d$)
in the proof of Lemma~\ref{simtemplemma} does not extend to the
boundary of $E_d$. (If it did, then the proof of
Theorem~\ref{simtempthm} would be complete simply by setting $h_d=h$
and applying Lemma~\ref{simtemplemma}.) However, the following lemma
shows that with the definition of $h_d$ used here, the extension to the
boundary does indeed hold.
% thus completing the proof of Theorem~\ref{simtempthm}.

%le9 #&#
\begin{lemma}\label{simtempboundarylemma}%
For all $h \in H$, for $x=1$ and for $x=\chi_d$,
\[
\lim_{d\to\infty} \bigl|d G^{(d)} h_d(x) - G^* h(x) \bigr|
= 0.
\]
\end{lemma}

\begin{pf}
We prove the case when $x=\chi_d$; the case $x=1$ is similar but
somewhat easier (since then $x$ does not depend on $d$).

Mimicking the Taylor expansion of Lemma~\ref{simtemplemma},
\begin{eqnarray*}
G^{(d)}h_d(\chi_d) &\stackrel{d}{\approx}&
{h_d'(\chi_d )[\alpha^-
(\chi_d^+ - \chi_d ) ]
\over2}\\
&&{} + {h_d''(\chi_d) \over4} \bigl[ \bigl(
\chi_d - \chi_d^+ \bigr) ^2 \alpha^- \bigr]
\\
&=& {h_d'(\chi_d) \over2}{\alpha^- \ell(\chi_d^+ )
% - \alpha^+ \ell(\beta_i^{(d)})
\over d^{1/2} } + {h_d''(\chi_d) \over4}
\biggl[ {\ell(\chi_d) ^2 \alpha^- \over d} \biggr]
\\
&\stackrel{d}{\approx}& {\alpha^- \ell(\chi_d^+) \over2 d^{1/2}} \bigl( h'(\chi)
+ \eta'(\chi) d^{-1/2} \bigr)\\
&&{} + {h_d''(\chi_d) \over4}
\biggl[ {\ell(\chi_d) ^2 \alpha^- \over d} \biggr].
\end{eqnarray*}
Thus since $h'(\chi) = 0$, this expression equals
\[
{h_d''(\chi_d) \over2} \biggl[ {\ell(\chi_d) ^2 \alpha^- \over d} \biggr].
\]
%
%$\lim_{d \to\infty} dG^{(d)} h (\chi_d ) - G^*h (\chi) =
%-h''(\chi) \ell^2(\chi) \lim_{d \to\infty} \alpha
%^-(\chi_d)/4$. %IS IT?? I AM WORRIED ABOUT THIS
%Now the above calculation shows that
%$$
%=
%$$
%
%Certainly $h_d \rightarrow h$ uniformly. However,

Next we note from~(\ref{alphaminus}) that
\[
\alpha^- \stackrel{1}{\approx} 2 \Phi \biggl( - {I^{1/2} \ell\over2}
\biggr).
\]
Hence, the above results show that
\[
\lim_{d \rightarrow\infty} d G_d h_d (
\chi_d) = \ell^2 (\chi) h''(
\chi) \Phi \biggl( - {I^{1/2} \ell\over2} \biggr).
\]
In light of formulae~(\ref{Geqn}) and~(\ref{sigmaeqn}), this completes
the proof.
\end{pf}

%The proof of Theorem~\ref{simtempthm} follows directly from
%Lemma~\ref{simtemplemma} %[RIGHT?].

Finally, we provide the missing proof from
Section~\ref{sec-mainproof}.

\begin{pf*}{Proof of Lemma~\ref{weakconvlemma}}
We first compute that, to first order as $h \searrow0$ and \mbox{$m \to
\infty$}, writing $x=i/m$ and $e=1/m$, we have
\begin{eqnarray*}
&& \mathbf{E} \biggl( Y^\sigma_{m, t + h} - Y^\sigma_{m, t } \Bigm|Y^\sigma_{m,t} = {i \over m} \biggr)
\\
&&\qquad \approx \biggl({m^2 S h \over2} \biggr) \biggl({1 \over m}
\biggr)\biggl({1 \over2S} \biggr)
\\
&&\qquad\quad{}\times \biggl[ \sigma^2
\biggl({i \over m} \biggr) + { \pi((i+1)/m)   \sigma^2 ((i+1)/m) \over\pi(i/m)}
\\
&&\hspace*{48pt}{} - \sigma^2 \biggl({i \over m} \biggr) -
{ \pi((i-1)/m)   \sigma^2 ((i-1)/m) \over\pi(i/m)} \biggr]
\\
&&\qquad = {h m \over4} \biggl[ { \pi(x+e)   \sigma^2(x+e) \over\pi(x) } -
{ \pi(x-e)   \sigma^2(x-e) \over\pi(x) } \biggr]
\\
%
% $$
% = % {h m \over4}
% \left[
% (\pi(x+e) \sigma^2(x+e)
% - \pi(x+e) \sigma^2(x-e))
% + (\pi(x+e) \sigma^2(x-e)
% - \pi(x-e) \sigma^2(x-e))
% \over\pi(x)
% \right]
% $$
%
&&\qquad\approx {h m \over4} \bigl[\bigl(\bigl(\pi(x)+e\pi'(x)
\bigr) \bigl(\sigma^2(x)+e \bigl(\sigma^2
\bigr)'(x) \bigr)
\\
&&\hspace*{56pt}{} - \bigl(\pi(x)-e\pi'(x) \bigr)
\bigl( \sigma^2(x)-e \bigl(\sigma^2 \bigr)'(x)
\bigr)\bigr)/\pi(x) \bigr]
\\
&&\qquad = {h m \over4} \biggl[{ 2e \pi'(x)
\sigma^2(x)+ 2e \pi(x) (\sigma^2)'(x) \over\pi(x)} \biggr]
\\
&&\qquad = {h m \over4} (2e) \bigl[ (\log\pi)'(x)
\sigma^2(x)+ 2 \sigma(x) \sigma'(x) \bigr]
\\
&&\qquad = h \biggl[ \frac{1}{2} (\log\pi)'(x) \sigma^2(x)+
\sigma(x) \sigma'(x) \biggr]
\end{eqnarray*}
and also
\begin{eqnarray*}
&& \mathbf{E} \biggl( \bigl(Y^\sigma_{m, t +
h} - Y^\sigma_{m, t } \bigr)^2 \Bigm| Y^\sigma_{m,t} = {i \over m}
\biggr)
\\
&&\qquad \approx \biggl({m^2 S h \over2} \biggr) \biggl( {1 \over2S}
\biggr) \biggl({1 \over m^2} \biggr) \bigl[2 \sigma^2(x) + 2
\sigma^2(x) \bigr] = h \bigl[ \sigma^2(x) \bigr].
\end{eqnarray*}
A comparison with~(\ref{first}) then shows that $Y^\sigma_m$ satisfies
the same first and second moment characteristics as $X^{\sigma}_t$, so
that $X^{\sigma}_t$ is indeed the correct putative limit.

% $$
% \approx % {h m \over4} \times
% \left[
% 2 \sigma^2({i \over m}) ((\log\pi)'({i \over m})
% + 4 \sigma({i \over m}) \sigma' ({i \over m}) )
% \right]
% $$
% $$
% = % {h m \over4} \times
% \left[
% 2 \sigma^2(x) ((\log\pi)'({i \over m})
% + 4 \sigma({i \over m}) \sigma' ({i \over m}) )
% \right]
% $$
% $$
% = % {2 \sigma^2({i \over m}) ((\log\pi)'({i \over m}) + 4 \sigma({i
% \over4 m }
% h
% $$
% $$
% = % {h \over2m} \sigma^2({i \over m}) ((\log\pi)'({i \over m})
% + {h \over m} \sigma({i \over m}) \sigma' ({i \over m}) )
%, % $$
% and also
% $$
% \mathbf{E} (
% (Y^\sigma_{m, t + h} - Y^\sigma_{m, t })^2 | Y^\sigma_{m,t} = i / m )
% \approx % {
% h m^2 S \times2 \sigma^2({i \over m}) \over
% 2 S m^2
% }
% = \sigma^2(i / m)   h
%  .
% $$
% Thus, the mean and variance of the process $Y^\sigma_m$ correspond to
% those of $\Xm$, to first order in $h$, as required.

In light of these calculations, the formal proof of this lemma then
proceeds along standard lines. Indeed, case~(a) is just a simpler
version of the proof of Theorem~\ref{simtempthm} above, and case~(b)
follows from standard arguments about using the uniform convergence of
generators (e.g.,~\cite{ethierkurtz}, Chapter~8) to establish the
approximation of birth and death processes by diffusions; see, for
example, Theorem~4.1 of Chapter~5 on page~387 of~\cite{bhat}.
\end{pf*}

%s4 #&#
\section{Discussion}

This paper has linked the usual Peskun ordering on asymptotic
variance of discrete-time Markov chains, to asymptotic variance of
diffusion processes. It has then applied these results to simulated
tempering algorithms, by proving that the inverse-temperatures of
such algorithms converge (in an appropriate limit) to a diffusion.
By maximising the speed of the resulting diffusion, it has obtained
results about the optimal choice of the temperature spacings.

We believe that Theorem~\ref{mainthm} could be useful in other contexts
as well, whenever we wish to compare two Langevin diffusion algorithms
directly, or alternatively whenever we wish to compare two
discrete-time processes which both have appropriate diffusion limits.

Of course, Theorem~\ref{mainthm} requires assumptions (A1) and (A2).
These are primarily just regularity assumptions, which would likely be
satisfied in most applications of interest. On the other hand, the
``exponentially-bounded tails'' aspect of assumption (A2) is more than
technical; rather, it provides us with some control over the extreme
tail excursions of the processes which we consider, and we suspect that
our limiting results might fail if no such control is provided.

Finally, our simulated tempering diffusion limit is only proven under
the rather strong and artificial assumption~(\ref{prodform}) involving
a product form of the target density. Indeed, this assumption is
central to our method of proof. However, as mentioned earlier, it is
known~\cite{RGG,statsci,bedard1,bedard2,bedardrosenthal,stuart} that
the general conclusions in this special case often hold in greater
generality, either approximately in numerical simulation studies, or
theoretically through more general methods of proof. In a similar
spirit, we believe that the simulated tempering diffusion limit proven
herein would approximately hold numerically in greater generality. In
addition, it might be possible to prove a stronger version of our
diffusion limit, with weaker assumptions, though such proofs would get
rather technical and we do not pursue them here.


%suskaldyti doi

% imsref loaded by linak, 2013-09-06 13:55:13

\printaddresses

\end{document}